\documentclass[a4paper, 12pt, twoside]{article}
\usepackage{latexsym,bm}
\usepackage[tbtags]{amsmath}
\usepackage{amssymb}

\usepackage{amsmath,amsfonts,amsthm,amssymb}

\usepackage{indentfirst}
\usepackage{paralist}
\usepackage[english]{babel}
\usepackage{amsmath,amssymb}

\usepackage{cite}
\usepackage{tabularx}
\usepackage{booktabs}
\usepackage{dcolumn}

\usepackage{multirow}

\usepackage{graphicx}

\usepackage{wrapfig}
\usepackage{multicol}
\usepackage{verbatim}

\newtheoremstyle{theorem}
  {10pt}
  {10pt}
  {\sl}
 {}
  {\bf}
  {. }
  { }
  {}
\theoremstyle{theorem}

\newtheorem{theorem}{Theorem}[section]

\newtheorem{definition}{Definition}[section]
 \newtheorem{lemma}{Lemma}[section]
 
 \newtheorem{remark}{Remark}[section]

\numberwithin{equation}{section}

\newtheoremstyle{defi}
  {10pt}
  {10pt}
  {\rm}
  {}
  {\bf}
  {. }
  { }
  {}
\theoremstyle{defi}

\allowdisplaybreaks \textheight 21 true cm
\usepackage{epsfig,graphicx,epsf,wrapfig}
\usepackage{dcolumn}
\usepackage{bm}
\begin{document} 


\title{{\large  \textbf{Reflected forward-backward stochastic differential equations driven by $G$-Brownian motion with continuous monotone coefficients \thanks{The work is supported in part by a NSFC Grant No. 11531006, NSFC Grant No. 11601203.}   }}
\author{Bingjun Wang$^{1,2}$
\footnote{Email: wbj586@126.com.cn.} \  \  Hongjun Gao$^{1}$ \footnote{Correspondence, Email: gaohj@njnu.edu.cn, gaohj@hotmail.com}\  \ Mei  Li$^{3}$ \footnote{ Email:  limei@njue.edu.cn
}
\\
\\
\footnotesize{ 1. Institute of  Mathematics, School of Mathematical Science} 
\\
\footnotesize{Nanjing Normal University, Nanjing 210023, P. R. China}
\\
\footnotesize{2. Jinling Institute of Technology, Nanjing 211169, P. R. China}
\\
\footnotesize{3. School of Applied Mathematics, Nanjing University of Finance and Economics}
\\
\footnotesize{ Nanjing 210023, P. R. China }}
 }  

 \date{}
\maketitle 
\footnotesize
\noindent \textbf{Abstract~~~} In this paper, we prove that there exists at least one solution for the  reflected forward-backward stochastic differential equation driven by $G$-Brownian motion satisfying the obstacle constraint with monotone coefficients.
\\[2mm]
\textbf{Key words~~~} reflected equation; forward-backward SDE; $G$-Brownian motion; monotone coefficients.
\\
\textbf{2010 Mathematics Subject Classification~~~}60H05, 60H10, 60H20

\section{Introduction}\label{s1}
Motivated by uncertainty problems, risk measures and the superhedging in finance, Peng systemically established a time-consistent fully expectation theory(see \cite{P5}). As a typical and important case, Peng introduced the $G$-expectation theory(see \cite{P1, P2}). In the $G$-expectation framework, the notion of $G$-Brownian motion and the corresponding stochastic calculus of It\^{o}'s type were established. On that basis, many properties and applications of the $G$-expectation, $G$-Brownian motion and the $G$-stochastic calculus are studied(see \cite{D1, G1}).

On that basis, some authors are interested in the forward stochastic differential equation driven by $G$-Brownian motion (FGSDE), which has a similar form as its counterpart in the classical framework, however, holds in a q.s. sense:
 $$X_{t}=x+\int_{0}^{t}b(s,X_{s})ds+\int_{0}^{t}h(s,X_{s})d\langle B\rangle_{s}+\int_{0}^{t}\sigma(s,X_{s})dB_{s},  \  0\leq t \leq T,\ q.s.,$$
where $\langle B\rangle$ is the quadratic variation of the $G$-Brownian motion $B$. Under the Lipschitz assumptions on the coefficients $b,h$ and $\sigma$, Peng \cite{P1} and Gao \cite{G1} have proved the wellposedness of such equation with the fixed-point iteration. Moreover, Bai and Lin \cite{B1} have studied the case when coefficients are integral-Lipschitz,  Lin \cite{L2} considered the reflected GSDEs with some good boundaries, Ren et al.\cite{R1} studied  stochastic functional differential equation with infinite delay driven by $G$-Brownian motion.

On the basis of a series of studies by Hu et al. \cite{H1} and Soner et al. \cite{S1} for the $G$-expectation, Peng et al. \cite{P3} obtained the complete representation theorem for $G$-martingale. Due to this contribution, Hu et al. \cite{H2}  obtained the existence, uniqueness, time consistency and a  priori estimates of  fully nonlinear backward stochastic differential equation  driven by a given $G$-Brownian motion (BGSDE) under standard Lipschitz conditions. Very recently, Li and Peng \cite{L4} study the reflected solution of  the following backward stochastic differential equations driven by $G$-Brownian motion (RBGSDE) via penalization:
  \begin{equation}
\begin{cases}
Y_{t}=\xi+\int_{t}^{T}f(s,X_{s},Y_{s}, Z_{s})ds+\int_{t}^{T}g(s,X_{s},Y_{s}, Z_{s})d\langle B\rangle _{s}-\int_{t}^{T}Z_{s}dB_{s}+(A_{T}-A_{t}),\\
Y_{t}\geq L_{t}, \  \  \{-\int_{0}^{t}(Y_{s}-L_{s})dA_{s}\}_{t\in [0,T]}     \   $is$ \   $a$   \   $non-increasing$  \    G-$martingale$.\nonumber
\end{cases}
\end{equation}
Under standard Lipschitz conditions on $f(s,y,z),g(s,y,z)$ in $y,z$ and the $L_{G}^{\beta}(\Omega_{T})(\beta>1)$ integrability condition on $\xi$,
there exists a triplet of processes $(Y,Z,A)\in \mathcal{S}_{G}^{\alpha}(0,T)$   satisfy above equation for $2\leq \alpha< \beta$. Here, $\mathcal{S}_{G}^{\alpha}(0,T)$ denote the collection of process $(Y,Z,A)$ such that $Y\in S_{G}^{\alpha}(0,T), Z\in H_{G}^{\alpha}(0,T)$, $A$ is a continuous nondecreasing process with $A_{0}=0$ and $A\in S_{G}^{\alpha}(0,T)$.

It is known that forward-backward equations are encountered when one applies the stochastic maximum principle to optimal stochastic control problems. Such equations are also encountered  in the probabilistic interpretation of a general type of systems quasilinear PDEs, as well as in finance (see \cite{H7} for example). In the linear expectation framework,  Antonelli et al. \cite{A1} and  Huang et al.\cite{H6} proved the existence of the solutions  for  backward-forward SDEs and reflected forward-backward SDEs respectively. However, in the $G$-framework, as far as  we know, there is no result about the reflected forward-backward stochastic differential equations driven by $G$-Brownian motion (RFBGSDEs) in which the solution of the BSDE stays above a given barrier.  One of the  differences is that the classical Skorohod condition should be  substituted by a $G$-martingale condition.  Moreover, the comparison theorem with respect to the increasing process $\{A_{t}\}_{t\in [0,T]}$ may not hold in  $G$-framework.  So  some  mathematical properties of this RFBGSDE should be developed.

In this paper,  we consider the solvability of the following RFBGSDEs with continuous monotone coefficients :
 \begin{equation}
\begin{cases}
X_{t}=x+\int_{0}^{t}b(s,X_{s},Y_{s})ds+\int_{0}^{t}h(s,X_{s},Y_{s})d\langle B\rangle _{s}+\int_{0}^{t}\sigma(s,X_{s})dB_{s},
\\
Y_{t}=\xi+\int_{t}^{T}f(s,X_{s},Y_{s}, Z_{s})ds+\int_{t}^{T}g(s,X_{s},Y_{s}, Z_{s})d\langle B\rangle _{s}-\int_{t}^{T}Z_{s}dB_{s}+(A_{T}-A_{t}),\\
Y_{t}\geq L_{t}, \  \  \{-\int_{0}^{t}(Y_{s}-L_{s})dA_{s}\}_{t\in [0,T]}     \   $is$ \   $a$   \   $non-increasing$  \    G-$martingale$.
\end{cases}
\end{equation}
 We  notice that the coefficients  of the forward GSDE contain the solution of the backward GSDE, so the forward GSDE  and the backward GSDE  are coupled together. Moreover, the coefficients only need to satisfy the linear growth condition, but  do not need to satisfy the Lipschitz condition.

The rest of this paper is organized as follows.  In section 2, we introduce some notions and results in the $G$-framework which are necessary for what follows.  In section 3, the existence theorem is provided.

\section{Preliminaries} \label{s2}
In this section, we introduce some notations and preliminary results in $G$-framework which are needed in the following sections. More details can be found in \cite{D1, G1, P1}.

Let $\Omega_{T}=C_{0}([0,T];R)$, the space of real valued continuous functions on $[0,T]$ with $w_{0}=0$, be endowed with the distance
$$d(w^{1},w^{2}):=\sum_{N=1}^{\infty}2^{-N}((\max_{0\leq t\leq N}|w_{t}^{1}-w_{t}^{2}|)\wedge1),$$
and let $B_{t}(w)=w_{t}$ be the canonical process. Denote by $\mathbb{F}:=\{\mathcal{F}_{t}\}_{0\leq t \leq T}$  the natural filtration generated by $B$, $L^{0}(\Omega_{T})$ be the space of all $\mathbb{F}$-measurable real functions.
 Let $L_{ip}(\Omega_{T}):=\{\varphi (B_{t_{1}},...,B_{t_{n}}): \forall n\geq1,t_{1},...,t_{n}\in [0,T], \forall\varphi\in C_{b,Lip}(R^{n})\},$  where $C_{b,Lip}(R^{n})$ denotes the set of bounded Lipschitz functions on $R^{n}$.   A sublinear functional on $L_{ip}(\Omega_{T})$  satisfies: for all $X,Y\in L_{ip}(\Omega_{T})$,

 (i) Monotonicity: $\mathbb{E}[X]\geq \mathbb{E}[Y]$ if $X\geq Y.$

 (ii) Constant preserving: $\mathbb{E}[C]=C$ for $C\in R.$

 (iii) Sub-additivity:  $\mathbb{E}[X+Y]\leq \mathbb{E}[X]+\mathbb{E}[Y].$

 (iv) Positive homogeneity: $\mathbb{E}[\lambda X]=\lambda \mathbb{E}[X]$ for  $\lambda\geq 0.$

The tripe $(\Omega,L_{ip}(\Omega_{T}),\mathbb{E})$ is called a sublinear expectation space and $\mathbb{E}$ is called a sublinear expectation.
\begin{definition}
A random variable $X\in L_{ip}(\Omega_{T})$ is $G$-normal distributed with parameters $(0, [\underline{\sigma}^{2}, \bar{\sigma}^{2}])$, i.e., $X\sim N(0, [\underline{\sigma}^{2}, \bar{\sigma}^{2}])$, if for each $\varphi\in C_{b,Lip}(R)$,  $u(t,x):=\mathbb{E}[\varphi(x+\sqrt{t}X)]$
is a viscosity solution to the following PDE on $R^{+}\times R$:
\begin{equation}
\begin{cases}
\frac{\partial u}{\partial t}+G(\frac{\partial^{2}u}{\partial x^{2}})=0,
\\
u_{t_{0}}=\varphi(x),
\end{cases}
\end{equation}
where $G(a):=\frac{1}{2}(a^{+}\bar{\sigma}^{2}-a^{-}{\underline{\sigma}}^{2}), a\in R$.
\end{definition}
\begin{definition}
We call a sublinear expectation $\hat{\mathbb{E}}:L_{ip}(\Omega_{T})\rightarrow R$ a $G$-expectation if the canonical process $B$ is a $G$-Brownian motion under $\hat{\mathbb{E}}[\cdot]$, that is, for each $0\leq s\leq t\leq T$, the increment $B_{t}-B_{s}\sim N(0, [\underline{\sigma}^{2}(t-s), \bar{\sigma}^{2}(t-s)])$ and for all $n>0, 0\leq t_{1}\leq \ldots \leq t_{n} \leq T$ and $\varphi\in L_{ip}(\Omega_{T})$,
$$\hat{\mathbb{E}}[\varphi(B_{t_{1}},\ldots,B_{t_{n-1}},B_{t_{n}}-B_{t_{n-1}} )]=\hat{\mathbb{E}}[\psi(B_{t_{1}},\ldots, B_{t_{n-1}})],$$
where $\psi(x_{1},\ldots,x_{n-1}):=\hat{\mathbb{E}}[\varphi(x_{1},\ldots,x_{n-1},\sqrt{t_{n}-t_{n-1}}B_{1})]$  and $B_{1}$ is $G$-normal distributed.
\end{definition}
For $p\geq 1$, we denote by $L_{G}^{p}(\Omega_{T})$ the completion of $L_{ip}(\Omega_{T})$ under the natural norm $\|X\|_{p,G}:=(\hat{\mathbb{E}}[|X|^{p}])^{\frac{1}{p}}$. $\hat{\mathbb{E}}$ is a continuous mapping on $L_{ip}(\Omega_{T})$ endowed with the norm $\|\cdot\|_{1,G}$. Therefore, it can be extended continuous to  $L_{G}^{1}(\Omega_{T})$ under the norm$\|X\|_{1,G}$.

Next, we introduce the It\^{o} integral of $G$-Brownian motion.

Let $M_{G}^{0}(0,T)$ be the collection of processes in the following form: for a given partition $\pi_{T}=\{t_{0},t_{1},...,t_{N}\}$ of $[0,T],$   set
$$\eta_{t}(w)=\sum_{k=0}^{N-1}\xi_{k}(w)I_{[t_{k},t_{k+1})}(t),$$
where $\xi_{k} \in L_{ip}(\Omega_{t_{k}}),k=0,1,...,N-1$ are given.   For $p\geq1$, we denote by  $H_{G}^{p}(0,T)$,    $M_{G}^{p}(0,T)$  the completion of  $M_{G}^{0}(0,T)$ under the norm
 $\|\eta\|_{H_{G}^{p}(0,T)}=\{\hat{\mathbb{E}}[(\int_{0}^{T}|\eta_{t}|^{2}dt)^{\frac{p}{2}}]\}^{\frac{1}{p}},$  $\|\eta\|_{M_{G}^{p}(0,T)}=\{\hat{\mathbb{E}}[\int_{0}^{T}|\eta_{t}|^{p}dt]\}^{\frac{1}{p}}$
 respectively.   It is easy to see that $H_{G}^{2}(0,T)=M_{G}^{2}(0,T)$.  Following Li and Peng \cite{L3}, for each $\eta\in H_{G}^{p}(0,T)$  with $p\geq 1$, we can define It\^{o} integral $\int_{0}^{T}\eta_{s}dB_{s}.$  Moreover, the following B-D-G inequality hold.
\begin{lemma}(\cite{G1})
Let $p\geq2$ and    $\eta\in M_{G}^{p}(0,T)$, then  we  have
$$\underline{\sigma}^{p}c_{p}\hat{\mathbb{E}}[(\int_{0}^{T}|\eta_{s}|^{2}ds)^{\frac{p}{2}}]\leq\hat{\mathbb{E}}[\sup_{0\leq t\leq T}|\int_{0}^{t}\eta_{s}dB_{s}|^{p}]\leq \overline{\sigma}^{p} C_{p}\hat{\mathbb{E}}[(\int_{0}^{T}|\eta_{s}|^{2}ds)^{\frac{p}{2}}],$$
where $0<c_{p}<C_{p}<\infty$ are constants.
\end{lemma}
Let $S_{G}^{0}(0,T)=\{h(t, B_{t_{1}\wedge t},\ldots,  B_{t_{n}\wedge t}): t, t_{1}, \ldots, t_{n} \in [0,T], h\in C_{b,lip}(R^{n+1})\}$.  For $p\geq1$  and $\eta\in S_{G}^{0}(0,T)$, set
$\|\eta\|_{S_{G}^{p}(0,T)}=(\hat{\mathbb{E}}[\sup_{0\leq t \leq T}|\eta_{t}|^{p}])^{\frac{1}{p}}$. Denote by $S_{G}^{p}(0,T)$ the completion of $S_{G}^{0}(0,T)$ under the norm $\|\eta\|_{S_{G}^{p}(0,T)}$.
\begin{definition}
Quadratic variation process of $G$-Brownian motion defined by
$$\langle B\rangle_{t}:=B_{t}^{2}-2\int_{0}^{t}B_{s}dB_{s}$$
is a continuous, nondecreasing process.
\end{definition}
For $\eta\in M_{G}^{0}(0,T)$, define $\int_{0}^{T}\eta_{s}d\langle B\rangle_{s}=\sum_{j=0}^{N-1}\xi_{j}(\langle B\rangle_{t_{j+1}}-\langle B\rangle_{t_{j}}): M_{G}^{0}(0,T)\rightarrow L_{G}^{1}(\Omega_{T})$. The mapping is continuous and can be extended to $M_{G}^{1}(0,T)$.
\begin{lemma} (\cite{P1, G1})
Let $p\geq 1$ and    $\eta\in M_{G}^{p}(0,T)$, then  we  have
$$\underline{\sigma}^{2}\hat{\mathbb{E}}[\int_{0}^{T}|\eta_{t}|dt]\leq \hat{\mathbb{E}}[|\int_{0}^{T}\eta_{t}d\langle B\rangle_{t}|]\leq\bar{\sigma}^{2}\hat{\mathbb{E}}[\int_{0}^{T}|\eta_{t}|dt],$$
$$\hat{\mathbb{E}}[\sup_{0\leq t\leq T}|\int_{0}^{t}\eta_{s}d\langle B\rangle _{s}|^{p}]\leq \overline{\sigma}^{p} C_{p}'\hat{\mathbb{E}}[\int_{0}^{T}|\eta_{s}|^{p}ds],$$
where $C_{p}'>0$ is a  constant independent of $\eta$.\end{lemma}

\begin{theorem}(\cite{D1})
There exists a weakly compact subset $\mathcal{P}\subset \mathcal{M}(\Omega_{T})$,  the set of probability measures on $(\Omega_{T},\mathcal{F}_{T})$, such that
$$\hat{\mathbb{E}}[\xi]=\max_{P\in\mathcal{ P}}E_{P}(\xi)  \  \ for  \  all \ \xi\in L_{G}^{1}(\Omega_{T}).$$
$\mathcal{ P}$ is called a set that represents $\hat{\mathbb{E}}$.
\end{theorem}
Let $\mathcal{ P}$ be a weakly compact set that represents $\hat{\mathbb{E}}$. For this $\mathcal{ P}$, we define capacity
$$c(A)=\sup_{P\in \mathcal{ P}}P(A), A\in \mathcal{F}_{T}.$$
A set $A\subset \Omega_{T}$ is a polar set  if $c(A)=0$. A property holds quasi-surely (q.s.) if it holds outside a polar set.

\begin{lemma}(\cite{D1})
Let $\{X^{n}\}_{n\in \mathbb{N}}\subset L_{G}^{1}(\Omega_{T})$ be such that $X^{n}\downarrow X$ q.s., then
$\hat{\mathbb{E}}[X^{n}]\downarrow \hat{\mathbb{E}}[X].$  In particular, if $X\in L_{G}^{1}(\Omega_{T})$, then $\hat{\mathbb{E}}[|X^{n}-X|]\downarrow 0 $, as $n\rightarrow\infty$.
 \end{lemma}

 \begin{lemma}(\cite{S2})
For any $\alpha \geq1$, $\delta>0$ and $1<\gamma <\beta:=\frac{(\alpha+\delta)}{\alpha}, \gamma \leq2$, we have
$$\hat{\mathbb{E}}[\sup_{t\in [0,T]}\hat{\mathbb{E}}_{t}[|\xi|^{\alpha}]]\leq \gamma^{*}\{(\hat{\mathbb{E}}|\xi|^{\alpha+\delta})^{\frac{\alpha}{\alpha+\delta}}+14^{\frac{1}{\gamma}}C_{\frac{\beta}{\gamma}}(\hat{\mathbb{E}}|\xi|^{\alpha+\delta})^{\frac{1}{\gamma}}\}, \     \   \forall \xi\in L_{ip}(\Omega_{T}),$$
where $C_{\frac{\beta}{\gamma}}=\sum_{i=1}^{\infty}i^{-\frac{\beta}{\gamma}},  \   \  \gamma^{*}=\frac{\gamma}{\gamma-1}.$
 \end{lemma}

To get the main result of this paper, we need the following  Lemma used by Lepeltier-San Martin.
\begin{lemma}(\cite{L5})
 Let $f:R^{m}\rightarrow R$ be a continuous function with linear growth, that is, there exist a constant $M<\infty$ such that $\forall x\in R^{m}, |f(x)|\leq M(1+|x|)$. Then the sequence of functions
 $$f_{n}(x)=\inf_{y\in Q^{p}}\{f(y)+n|x-y|\}$$
 is well defined for $n\geq M$ and satisfies

 (i) linear growth: $\forall x\in R^{m}, |f_{n}(x)|\leq M(1+|x|)$;

 (ii) monotonicity in $n$: $\forall x\in R^{m}, f_{n}(x)\leq f_{n+1}(x)$;

 (iii) Lipschitz condition:  $\forall x, y \in R^{m}, |f_{n}(x)-f_{n}(y)|\leq n|x-y|$;

 (iv) strong convergence: if $x_{n}\rightarrow x$, then  $f_{n}(x_{n})\rightarrow f(x)$.
 \end{lemma}

\section{main result}  \label{s3}
 \begin{definition}
 A quadruple of processes $(X,Y,Z,A)$  is called a solution of  reflected  FBGSDEs (1.1) if the following properties are satisfied:

 (i) $X\in M_{G}^{2}(0,T)$, $(Y,Z,A)\in \mathcal{S}_{G}^{2}(0,T)$;

 (ii)  \begin{equation}
\begin{cases}
X_{t}=x+\int_{0}^{t}b(s,X_{s},Y_{s})ds+\int_{0}^{t}h(s,X_{s},Y_{s})d\langle B\rangle _{s}+\int_{0}^{t}\sigma(s,X_{s})dB_{s},
\\
Y_{t}=\xi+\int_{t}^{T}f(s,X_{s},Y_{s}, Z_{s})ds+\int_{t}^{T}g(s,X_{s},Y_{s}, Z_{s})d\langle B\rangle _{s}-\int_{t}^{T}Z_{s}dB_{s}+(A_{T}-A_{t});\nonumber
\end{cases}
\end{equation}

 (iii)  $Y_{t}\geq L_{t}$, \  \  $\{-\int_{0}^{t}(Y_{s}-L_{s})dA_{s}\}_{t\in [0,T]}$      is    a     non-increasing     $G$-martingale.
 \end{definition}
In the sequel, we will work under the following assumptions: for any $s\in [0,T], w\in\Omega, x,x',y,z\in R,  \beta>2$:

(H1) $b(\cdot,x,y,z),h(\cdot,x,y,z), \sigma(\cdot, x)\in M_{G}^{2}(0,T)$, $f(\cdot,x,y,z),g(\cdot,x,y,z)\in M_{G}^{\beta}(0,T)$;

(H2)  $b,h$ are incresing in $y$ and $f,g$ are incresing in $x$;

(H3)  there exists a constant $M>0$, such that
$$|b(s,x,y)|\vee |h(s,x,y)|\leq M(1+|x|+|y|),  \  \   |f(s,x,y,z)|\vee |g(s,x,y,z)|\leq M(1+|y|+|z|);$$
$$|\sigma(s,x)|\leq M(1+|x|),  \  \   |\sigma(s,x)-\sigma(s,x')|\leq M|x-x'|;$$

(H4)   $\xi\in L_{G}^{\beta}(\Omega_{T})$ and $\xi\geq L_{T}$, q.s..

(H5) $\{L_{t}\}_{t\in[0,T]}\in S_{G}^{\beta}(0,T)$ and there exists a constant $c$ such that  $L_{t}\leq c$, for each $t\in[0,T]$;

 For notational simplification, by Lemma 2.2, we only consider the case $h=0$ and $g=0$. I.e., we consider the following equation:
  \begin{equation}
\begin{cases}
X_{t}=x+\int_{0}^{t}b(s,X_{s},Y_{s})ds+\int_{0}^{t}\sigma(s,X_{s})dB_{s},
\\
Y_{t}=\xi+\int_{t}^{T}f(s,X_{s},Y_{s}, Z_{s})ds-\int_{t}^{T}Z_{s}dB_{s}+(A_{T}-A_{t}),\\
Y_{t}\geq L_{t}, \  \  \{-\int_{0}^{t}(Y_{s}-L_{s})dA_{s}\}_{t\in [0,T]}     \   $is$ \   $a$   \   $non-increasing$  \    G-$martingale$.
\end{cases}
\end{equation}
But the results still hold for other case.  In the following, $C$ always denote a positive constant which  may change from line to line.
 \begin{theorem}
Suppose that $\xi, b,f, \sigma$ satisfy (H1)-(H4), $L$ satisfies (H5). Then the RFBGSDE (3.1) has at least one solution $(X,Y,Z,A)$.
\end{theorem}
 {\bf Proof~~~}  In order to construct a solution of (3.1), our basic idea is to consider the following iteration:
  \begin{equation}
\begin{cases}
X^{n}_{t}=x+\int_{0}^{t}b(s,X^{n}_{s},Y^{n}_{s})ds+\int_{0}^{t}\sigma(s,X^{n}_{s})dB_{s},
\\
Y^{n}_{t}=\xi+\int_{t}^{T}f(s,X^{n-1}_{s},Y^{n}_{s},  Z^{n}_{s})ds-\int_{t}^{T}Z^{n}_{s}dB_{s}+(A^{n}_{T}-A^{n}_{t}),\\
Y^{n}_{t}\geq L_{t}, \{-\int_{0}^{t}(Y^{n}_{s}-L_{s})dA^{n}_{s}\}_{t\in [0,T]}\      $is$ \      $a$ \      $non-increasing$    \  G-$martingale$.
\end{cases}
\end{equation}
We will show that the limit of the sequence $\{(X^{n},Y^{n},Z^{n},A^{n})\}_{n\in \mathbb{N}}$  verifies equations  (3.1).

Step 1: Construction of the starting point.

Let us consider the following two standard  reflected backward $G$-stochastic differential equations:
  \begin{equation}
\begin{cases}
Y_{t}^{0}=\xi-K\int_{t}^{T}(1+|Y_{s}^{0}|+|Z_{s}^{0}|)ds-\int_{t}^{T}Z_{s}^{0}dB_{s}+(A_{T}^{0}-A_{t}^{0}), \\
Y^{0}_{t}\geq L_{t}, \{-\int_{0}^{t}(Y_{s}^{0}-L_{s})dA_{s}^{0}\}_{t\in [0,T]}\      $is$ \      $a$ \      $non-increasing$    \  G-$martingale$
\end{cases}
\end{equation}
and
  \begin{equation}
\begin{cases}
U_{t}=|\xi|+K\int_{t}^{T}(1+|U_{s}|+|V_{s}|)ds-\int_{t}^{T}V_{s}dB_{s}+(N_{T}-N_{t}),\\
U_{t}\geq L_{t}, \{-\int_{0}^{t}(U_{s}-L_{s})dN_{s}\}_{t\in [0,T]}\     $is$ \      $a$ \      $non-increasing$    \  G-$martingale$,\nonumber
\end{cases}
\end{equation}
where $K>0$  is a constant. By virtue of the Lipschitz property of the coefficients $\pm K(1+|y|+|z|)$, thanks to Theorem 5.1 in \cite{L4}, each one has an unique solution denoted by $(Y^{0},Z^{0},A^{0})$ and $(U,V,N)$ respectively. More precisely,  $(Y^{0}_{t}, Z^{0}_{t}, A_{t}^{0}), (U_{t}, V_{t}, N_{t}) \in \mathcal{S}_{G}^{\alpha}(0,T)$ for  $2\leq \alpha < \beta$,   such that $
\{-\int_{0}^{t}(Y_{s}^{0}-L_{s})dA_{s}^{0}\}_{t\in [0,T]}, \{-\int_{0}^{t}(U_{s}-L_{s})dN_{s}\}_{t\in [0,T]}$  are  non-increasing $G$- martingale.  By the  comparison theorem in \cite{L4}, we know that for all $t\in [0,T]$, $Y_{t}^{0}\leq U_{t}$,  q.s..

Step 2: Construction of $X^{0}$.

Now, we consider the forward equation
\begin{align}
&X_{t}^{0}=x+\int_{0}^{t}b(s,X^{0}_{s},Y_{s}^{0})ds+\int_{0}^{t}\sigma(s,X^{0}_{s})dB_{s},
\end{align}
where $Y^{0}$ is the solution of (3.3).

Let $\{b_{k}(s,x,y)\}_{k\geq 0}$  be the sequence defined in Lemma 2.5.  then we can conclude that the following GSDE has a unique solution $X_{t}^{0,k}\in M_{G}^{2}(0,T)$  by the lipschitz property of $b_{k}$, i.e.,
\begin{align}
&X_{t}^{0,k}=x+\int_{0}^{t}b_{k}(s,X^{0,k}_{s},Y_{s}^{0})ds+\int_{0}^{t}\sigma(s,X^{0,k}_{s})dB_{s}.
\end{align}
Moreover, by the comparison theorem in  \cite{L1} and Lemma 2.5, we know that for $t\in [0,T]$,  $X_{t}^{0,k}\leq X_{t}^{0,k+1}\leq S_{t}$,  where $S_{t}\in M_{G}^{2}(0,T)$ is the unique solution of the following GSDE:
\begin{align}
&S_{t}=x+K\int_{0}^{t}(1+|S_{s}|+|U_{s}|)ds+\int_{0}^{t}\sigma(s,S_{s})dB_{s}.\nonumber
\end{align}
Actually,  by  Corollary 3.2 in chapter V of  \cite{P1}, we know that  $X_{t}^{0,k}, S_{t}\in S_{G}^{2}(0,T)$.  So there exists a lower semi-continuous process $X_{t}^{0}\in S_{G}^{2}(0,T)$ such that $X_{t}^{0,k}\uparrow X_{t}^{0}$  as $k\rightarrow\infty$,   q.s.. Notice that $X_{t}^{0,k}, X_{t}^{0}\in S_{G}^{2}(0,T)$, which obviously belong to a larger space $L_{G}^{2}(0,T)$. Then by the downward   monotone convergence theorem (Lemma 2.3), we have $\hat{\mathbb{E}}[|X_{t}^{0,k}- X_{t}^{0}|^{2}]\downarrow 0$ as $k\rightarrow\infty$.

By  Lemma 2.5 and the dominated convergence theorem with respect to $t$, we have
\begin{align}
&\hat{\mathbb{E}}[\int_{0}^{T}|b_{k}(s,X_{s}^{0,k},Y_{s}^{0})-b(s,X_{s}^{0},Y_{s}^{0})|^{2}ds]\nonumber\\
&\leq\hat{\mathbb{E}}[\int_{0}^{T}|b_{k}(s,X_{s}^{0,k},Y_{s}^{0})-b_{k}(s,X_{s}^{0},Y_{s}^{0})|^{2}ds]+\hat{\mathbb{E}}[\int_{0}^{T}|b_{k}(s,X_{s}^{0},Y_{s}^{0})-b(s,X_{s}^{0},Y_{s}^{0})|^{2}ds]\nonumber\\
&\leq C\int_{0}^{T}\hat{\mathbb{E}}[|X_{s}^{0,k}-X_{s}^{0}|^{2}]ds+\int_{0}^{T}\hat{\mathbb{E}}[|b_{k}(s,X_{s}^{0},Y_{s}^{0})-b(s,X_{s}^{0},Y_{s}^{0})|^{2}]ds\rightarrow0,
\end{align}
as $k\rightarrow\infty$.

On the other hand, since $|\sigma(s.X_{s}^{0,k})-\sigma(s.X_{s}^{0})|\leq M|X_{s}^{0,k}-X_{s}^{0}|$, then according to Lemma 2.1, we have
\begin{align}
&\hat{\mathbb{E}}[\sup_{0\leq t \leq T}\int_{0}^{t}|\sigma(s,X_{s}^{0,k})-\sigma(s,X_{s}^{0})|^{2}ds]\leq C\int_{0}^{T}\hat{\mathbb{E}}|X_{s}^{0,k}-X_{s}^{0}|^{2}ds\rightarrow0,\nonumber
\end{align}
as $k\rightarrow\infty$.

Now taking limit on both side of (3.5), then we obtain that the continuous process  $\{X_{t}^{0}\}_{t\in [0,T]}$ satisfies (3.4).

Step 3: Construction of $( X^{n}, Y^{n}, Z^{n}, A^{n})$.

We  focus on $(X^{1}, Y^{1}, Z^{1}, A^{1})$. First,  based on $X^{0}$, we can construct $Y^{1}$. In fact, denote $f^{1}(s,w,y,z)=f(s,X_{s}^{0}(w), y,z)$, then  one can easily check that $|f^{1}(s,w,y,z)|\leq K(1+|y|+|z|)$. Define once again $f^{1}_{k}(s,w,y,z)$ the approximating sequence in Lemma 2.5, then by  Theorem 5.1 in \cite{L4}, for $2\leq \alpha <\beta$, we have a unique triple $(Y^{1,k}_{\cdot}, Z^{1,k}_{\cdot}, A^{1,k}_{\cdot})\in \mathcal{S}_{G}^{\alpha}(0,T)$ satisfying
  \begin{equation}
\begin{cases}
Y_{t}^{1,k}=\xi+\int_{t}^{T}f^{1}_{k}(s,w,Y_{s}^{1,k},Z_{s}^{1,k})ds-\int_{t}^{T}Z_{s}^{1,k}dB_{s}+(A_{T}^{1,k}-A_{t}^{1,k}), \\
Y^{1,k}_{t}\geq L_{t}, \{-\int_{0}^{t}(Y_{s}^{1,k}-L_{s})dA_{s}^{1,k}\}_{t\in [0,T]}\      $is$ \      $a$ \      $non-increasing$    \  G-$martingale$.
\end{cases}
\end{equation}
Moreover,  by the  comparison theorem in \cite{L4}, we have
$$Y_{t}^{0}\leq Y_{t}^{1,k}\leq Y_{t}^{1,k+1}\leq U_{t},   \  \  t\in [0,T],  \   q.s..$$
Then it is easy to see that  there exists a constant $C$ independent of  $k$, such that
\begin{align}
&\hat{\mathbb{E}} [\sup_{0\leq t\leq T} |Y_{t}^{1,k}|^{\alpha}]\leq C.\end{align}

By Proposition 3.1 in \cite{L4}, we have
\begin{align}
&\hat{\mathbb{E}} [(\int_{0}^{T}|Z_{s}^{1,k}|^{2}ds)^{\frac{\alpha}{2}}]\leq C'\{\hat{\mathbb{E}}[\sup_{0\leq t\leq T}|Y_{t}^{1,k}|^{\alpha}]+M^{\frac{\alpha}{2}}(\hat{\mathbb{E}}[\sup_{0\leq t\leq T} |Y_{t}^{1,k}|^{\alpha}])^{\frac{1}{2}}\}\leq C
\end{align}
and
\begin{align}
&\hat{\mathbb{E}} [|A_{T}^{1,k}|^{\alpha}]\leq C_{\alpha}\{\hat{\mathbb{E}}\sup_{0\leq t\leq T}|Y_{t}^{1,k}|^{\alpha}+(MT)^{\alpha}\}\leq C,
\end{align}
where $C'$ is a constant.

 As $\{Y^{1,k}\}_{k\geq1}$  is an increasing sequence, we denote the limit by $Y^{1}$.  It is easy to see that $Y_{t}^{1}\geq Y_{t}^{0}$ for each $t\in [0,T]$.   Moreover, by Fatou's lemma, we have $Y_{t}^{1}\in S_{G}^{\alpha}(0,T)$. Hence, by the dominated convergence theorem with respect to $t$, we have
 $$\int_{0}^{T}\hat{\mathbb{E}}|Y^{1,k}_{s}-Y_{s}^{1}|^{\alpha}ds\rightarrow 0,$$
 as $k\rightarrow\infty$.  i.e., $Y^{1,k}\rightarrow Y^{1}$ in $M_{G}^{\alpha}(0,T)$. In the following, we can show that this convergence holds in $S_{G}^{\alpha}(0,T)$.

 Applying $G$-It\^{o}'s formula to $(|Y^{1,k}_{t}-Y^{1,j}_{t}|^{2})^{\frac{\alpha}{2}}$, we know that
 \begin{align}
&|Y^{1,k}_{t}-Y^{1,j}_{t}|^{\alpha}+\frac{\alpha}{2}\int_{t}^{T}(|Y^{1,k}_{t}-Y^{1,j}_{t}|^{2})^{\frac{\alpha}{2}-1}|Z^{1,k}_{s}-Z^{1,j}_{s}|^{2}d\langle B\rangle_{s}\nonumber\\
&=\alpha(1-\frac{\alpha}{2})\int_{t}^{T}(|Y^{1,k}_{t}-Y^{1,j}_{t}|^{2})^{\frac{\alpha}{2}-2}(Y^{1,k}_{t}-Y^{1,j}_{t})^{2}|Z^{1,k}_{s}-Z^{1,j}_{s}|^{2}d\langle B\rangle_{s}\nonumber\\
&+\alpha\int_{t}^{T}(|Y^{1,k}_{s}-Y^{1,j}_{s}|^{2})^{\frac{\alpha}{2}-1}(Y^{1,k}_{s}-Y^{1,j}_{s})[f_{k}^{1}(s,w,Y_{s}^{1,k},Z_{s}^{1,k})-f_{j}^{1}(s,w,Y_{s}^{1,j},Z_{s}^{1,j})]ds\nonumber\\
&+\alpha\int_{t}^{T}(|Y^{1,k}_{s}-Y^{1,j}_{s}|^{2})^{\frac{\alpha}{2}-1}(Y^{1,k}_{s}-Y^{1,j}_{s})d(A^{1,k}_{s}-A^{1,j}_{s})\nonumber\\
&-\alpha\int_{t}^{T}(|Y^{1,k}_{s}-Y^{1,j}_{s}|^{2})^{\frac{\alpha}{2}-1}(Y^{1,k}_{s}-Y^{1,j}_{s})(Z^{1,k}_{s}-Z^{1,j}_{s})dB_{s}\nonumber\\
&\leq\alpha\int_{t}^{T}|Y^{1,k}_{s}-Y^{1,j}_{s}|^{\alpha-1}|f_{k}^{1}(s,w,Y_{s}^{1,k},Z_{s}^{1,k})-f_{j}^{1}(s,w,Y_{s}^{1,j},Z_{s}^{1,j})|ds\nonumber\\
&+\alpha\int_{t}^{T}(|Y^{1,k}_{s}-Y^{1,j}_{s}|^{2})^{\frac{\alpha}{2}-1}(Y^{1,k}_{s}-Y^{1,j}_{s})^{+}dA^{1,k}_{s}\nonumber\\
&+\alpha\int_{t}^{T}(|Y^{1,k}_{s}-Y^{1,j}_{s}|^{2})^{\frac{\alpha}{2}-1}(Y^{1,k}_{s}-Y^{1,j}_{s})^{-}dA^{1,j}_{s}\nonumber\\
&-\alpha\int_{t}^{T}(|Y^{1,k}_{s}-Y^{1,j}_{s}|^{2})^{\frac{\alpha}{2}-1}(Y^{1,k}_{s}-Y^{1,j}_{s}) (Z^{1,k}_{s}-Z^{1,j}_{s})dB_{s}.
\end{align}

Let $M^{j,k}_{t}=-\alpha\int_{0}^{t}(|Y^{1,k}_{s}-Y^{1,j}_{s}|^{2})^{\frac{\alpha}{2}-1}(Y^{1,k}_{s}-Y^{1,j}_{s})^{+}dA^{1,k}_{s}-\alpha\int_{0}^{t}(|Y^{1,k}_{s}-Y^{1,j}_{s}|^{2})^{\frac{\alpha}{2}-1}(Y^{1,k}_{s}-Y^{1,j}_{s})^{-}dA^{1,j}_{s}+\alpha\int_{0}^{t}(|Y^{1,k}_{s}-Y^{1,j}_{s}|^{2})^{\frac{\alpha}{2}-1}(Y^{1,k}_{s}-Y^{1,j}_{s}) (Z^{1,k}_{s}-Z^{1,j}_{s})dB_{s}$.
Since
$$0\geq-(Y^{1,k}_{s}-Y^{1,j}_{s})^{+}\geq-(Y^{1,k}_{s}-L_{s}),$$
then
$$0\geq\int_{0}^{t}-(Y^{1,k}_{s}-Y^{1,j}_{s})^{+}dA^{1,k}_{s}\geq-\int_{0}^{t}(Y^{1,k}_{s}-L_{s})A^{1,k}_{s},$$
which implies $\{M_{t}^{j,k}\}_{t\in[0,T]}$  is a $G$-martingale. We rewrite (3.11) as
 \begin{align}
&M^{j,k}_{T}-M^{j,k}_{t}+|Y^{1,k}_{t}-Y^{1,j}_{t}|^{\alpha}+\frac{\alpha}{2}\int_{t}^{T}|Y^{1,k}_{t}-Y^{1,j}_{t}|^{\alpha-2}|Z^{1,k}_{s}-Z^{1,j}_{s}|^{2}d\langle B\rangle_{s}\nonumber\\
&\leq\alpha\int_{t}^{T}|Y^{1,k}_{s}-Y^{1,j}_{s}|^{\alpha-1}|f_{k}^{1}(s,w,Y_{s}^{1,k},Z_{s}^{1,k})-f_{j}^{1}(s,w,Y_{s}^{1,j},Z_{s}^{1,j})|ds.\nonumber
\end{align}
Taking conditional expectation on both side, we obtain
\begin{align}
&|Y^{1,k}_{t}-Y^{1,j}_{t}|^{\alpha}+\frac{\alpha}{2}\hat{\mathbb{E}}_{t}[\int_{t}^{T}|Y^{1,k}_{t}-Y^{1,j}_{t}|^{\alpha-2}|Z^{1,k}_{s}-Z^{1,j}_{s}|^{2}d\langle B\rangle_{s}]\nonumber\\
&\leq\alpha\hat{\mathbb{E}}_{t}[\int_{t}^{T}|Y^{1,k}_{s}-Y^{1,j}_{s}|^{\alpha-1}|f_{k}^{1}(s,w,Y_{s}^{1,k},Z_{s}^{1,k})-f_{j}^{1}(s,w,Y_{s}^{1,j},Z_{s}^{1,j})|ds].
\end{align}
By Lemma 2.4, for $0<\varepsilon<\frac{4-\alpha}{3\alpha-4},$  we have
\begin{align}
&\hat{\mathbb{E}}[\sup_{0\leq t\leq T}|Y^{1,k}_{t}-Y^{1,j}_{t}|^{\alpha}]\nonumber\\
&\leq\alpha\hat{\mathbb{E}}[\int_{0}^{T}|Y^{1,k}_{s}-Y^{1,j}_{s}|^{\frac{\alpha}{2}}|Y^{1,k}_{s}-Y^{1,j}_{s}|^{\frac{\alpha}{2}-1}|f_{k}^{1}(s,w,Y_{s}^{1,k},Z_{s}^{1,k})-f_{j}^{1}(s,w,Y_{s}^{1,j},Z_{s}^{1,j})|ds]^{1+\varepsilon}\nonumber\\
&\leq\alpha(\hat{\mathbb{E}}\int_{0}^{T}|Y^{1,k}_{s}-Y^{1,j}_{s}|^{\alpha}ds)^{\frac{1+\varepsilon}{2}}\nonumber\\
&\cdot[\hat{\mathbb{E}}\sup_{t\in[0,T]}|Y^{1,k}_{t}-Y^{1,j}_{t}|^{\frac{\alpha}{2}}(\int_{0}^{T}|f_{k}^{1}(s,w,Y_{s}^{1,k},Z_{s}^{1,k})-f_{j}^{1}(s,w,Y_{s}^{1,j},Z_{s}^{1,j})|^{2}ds)^{\frac{1+\varepsilon}{1-\varepsilon}}]^{\frac{1-\varepsilon}{2}}\nonumber\\
&\leq\alpha(\hat{\mathbb{E}}\int_{0}^{T}|Y^{1,k}_{s}-Y^{1,j}_{s}|^{\alpha}ds)^{\frac{1+\varepsilon}{2}}\nonumber\\
&\cdot\{(\hat{\mathbb{E}}[\sup_{t\in[0,T]}|Y^{1,k}_{t}-Y^{1,j}_{t}|^{\alpha}])^{\frac{1}{2}}(\hat{\mathbb{E}}[\int_{0}^{T}|f_{k}^{1}(s,w,Y_{s}^{1,k},Z_{s}^{1,k})-f_{j}^{1}(s,w,Y_{s}^{1,j},Z_{s}^{1,j})|^{2}ds]^{\frac{2(1+\varepsilon)}{1-\varepsilon}})^{\frac{1}{2}}\}^{\frac{1-\varepsilon}{2}}.\nonumber
\end{align}
By (3.8)-(3.10) and Lemma 2.5, it is easy to see that the second term in the above formula  is bounded. Then we have
\begin{align}
&\hat{\mathbb{E}}[\sup_{0\leq t\leq T}|Y^{1,k}_{t}-Y^{1,j}_{t}|^{\alpha}]\leq C(\hat{\mathbb{E}}\int_{0}^{T}|Y^{1,k}_{s}-Y^{1,j}_{s}|^{\alpha}ds)^{\frac{1+\varepsilon}{2}}\rightarrow0,
\end{align}
as $j,k\rightarrow\infty$. Thus $\{Y^{1,k}_{t}\}_{k\in \mathbb{N}}$ is a Cauchay sequence in $S_{G}^{\alpha}(0,T)$ with the limit $Y_{t}^{1}$.

Let  $\alpha=2$ in (3.11), by Lemma 2.2,  we have
 \begin{align}
&\underline{\sigma}^{2}\int_{0}^{T}|Z^{1,k}_{s}-Z^{1,j}_{s}|^{2}ds\leq C\{\int_{0}^{T}|Y^{1,k}_{s}-Y^{1,j}_{s}||f_{k}^{1}(s,w,Y_{s}^{1,k},Z_{s}^{1,k})-f_{j}^{1}(s,w,Y_{s}^{1,j},Z_{s}^{1,j})|ds\nonumber\\
&+\int_{0}^{T}|Y^{1,k}_{s}-Y^{1,j}_{s}|d(A_{s}^{1,j}+A_{s}^{1,k})-\int_{0}^{T}(Y^{1,k}_{s}-Y^{1,j}_{s})(Z^{1,k}_{s}-Z^{1,j}_{s})dB_{s}\}\nonumber\\
&\leq C\{\sup_{t\in [0,T]} |Y^{1,k}_{s}-Y^{1,j}_{s}|\cdot\int_{0}^{T}|f_{k}^{1}(s,w,Y_{s}^{1,k},Z_{s}^{1,k})-f_{j}^{1}(s,w,Y_{s}^{1,j},Z_{s}^{1,j})|ds\nonumber\\
&+\sup_{t\in [0,T]} |Y^{1,k}_{s}-Y^{1,j}_{s}|\cdot(|A_{T}^{1,j}|+|A_{T}^{1,k}|)-\int_{0}^{T}(Y^{1,k}_{s}-Y^{1,j}_{s})(Z^{1,k}_{s}-Z^{1,j}_{s})dB_{s}\},
\end{align}
 By Lemma 2.1, for any $\varepsilon>0$, we obtain
  \begin{align}
&\hat{\mathbb{E}}[(\int_{0}^{T}(Y^{1,k}_{s}-Y^{1,j}_{s})(Z^{1,k}_{s}-Z^{1,j}_{s})dB_{s})^{\frac{\alpha}{2}}]\leq C\hat{\mathbb{E}}[(\int_{0}^{T}(Y^{1,k}_{s}-Y^{1,j}_{s})^{2}(Z^{1,k}_{s}-Z^{1,j}_{s})^{2}ds)^{\frac{\alpha}{4}}]\nonumber\\
&\leq C(\hat{\mathbb{E}}[\sup_{t\in [0,T]} |Y^{1,k}_{s}-Y^{1,j}_{s}|^{\alpha}])^{\frac{1}{2}}(\hat{\mathbb{E}}[(\int_{0}^{T} |Z^{1,k}_{s}-Z^{1,j}_{s}|^{2}ds)^{\frac{\alpha}{2}}])^{\frac{1}{2}}\nonumber\\
&\leq \frac{C}{4\varepsilon}\hat{\mathbb{E}}[\sup_{t\in [0,T]} |Y^{1,k}_{s}-Y^{1,j}_{s}|^{\alpha}]+C\varepsilon\hat{\mathbb{E}}[(\int_{0}^{T} |Z^{1,k}_{s}-Z^{1,j}_{s}|^{2}ds)^{\frac{\alpha}{2}}].
\end{align}
By (3.8)-(3.10) and the H\"{o}lder inequality, choosing $\varepsilon$ small enough, it follows from (3.14) that
\begin{align}
&\hat{\mathbb{E}}[(\int_{0}^{T}|Z^{1,k}_{s}-Z^{1,j}_{s}|^{2}ds)^{\frac{\alpha}{2}}]\nonumber\\
&\leq C\{(\hat{\mathbb{E}}[\sup_{t\in [0,T]} |Y^{1,k}_{s}-Y^{1,j}_{s}|^{\alpha}])^{\frac{1}{2}}(\hat{\mathbb{E}}[(\int_{0}^{T}|f_{k}^{1}(s,w,Y_{s}^{1,k},Z_{s}^{1,k})-f_{j}^{1}(s,w,Y_{s}^{1,j},Z_{s}^{1,j})|^{2}ds)^{\frac{\alpha}{2}}])^{\frac{1}{2}}\nonumber\\
&+(\hat{\mathbb{E}}[\sup_{t\in [0,T]} |Y^{1,k}_{s}-Y^{1,j}_{s}|^{\alpha}])^{\frac{1}{2}}(\hat{\mathbb{E}}|A_{T}^{1,k}|^{\alpha}+\hat{\mathbb{E}}|A_{T}^{1,j}|^{\alpha})^{\frac{1}{2}}+\hat{\mathbb{E}}[\sup_{t\in [0,T]} |Y^{1,k}_{s}-Y^{1,j}_{s}|^{\alpha}]\}
\nonumber\\
&\leq C\{\hat{\mathbb{E}}[\sup_{t\in [0,T]} |Y^{1,k}_{s}-Y^{1,j}_{s}|^{\alpha}]+(\hat{\mathbb{E}}[\sup_{t\in [0,T]} |Y^{1,k}_{s}-Y^{1,j}_{s}|^{\alpha}])^{\frac{1}{2}}\}.\nonumber
\end{align}
It is straightforward to   show that
\begin{align}
&\lim_{j,k\rightarrow\infty}\hat{\mathbb{E}}[(\int_{0}^{T}|Z^{1,k}_{s}-Z^{1,j}_{s}|^{2}ds)^{\frac{\alpha}{2}}]=0.
\end{align}
Then there exists a process $Z_{t}^{1}\in H_{G}^{\alpha}(0,T)$ such that
$\hat{\mathbb{E}}[(\int_{0}^{T}|Z^{1,k}_{s}-Z^{1}_{s}|^{2}ds)^{\frac{\alpha}{2}}]\rightarrow 0,$
as $k \rightarrow\infty$.

Similar  to that in (3.6), it is easy to see that
 \begin{align}
&\hat{\mathbb{E}}[|\int_{0}^{T}f^{1}_{k}(s,w,Y_{s}^{1,k},Z_{s}^{1,k})-f^{1}(s,w,Y_{s}^{1},Z_{s}^{1})ds|^{\alpha}]\rightarrow0,   \       \    \   as \   \  k\rightarrow\infty.\end{align}

Meanwhile, since
 \begin{align}
&A_{t}^{1,k}-A_{t}^{1,j}=Y_{0}^{1,k}-Y_{0}^{1,j}-(Y_{t}^{1,k}-Y_{t}^{1,j})\nonumber\\
&-\int_{0}^{t}[f_{k}^{1}(s,w,Y_{s}^{1,k},Z_{s}^{1,k})-f_{j}^{1}(s,w,Y_{s}^{1,j},Z_{s}^{1,j})]ds+\int_{0}^{t}(Z_{s}^{1,k}-Z_{s}^{1,j})dB_{s}.\nonumber
\end{align}
Then by Lemma 2.1, (3.13), (3.16) and (3.17), we have
\begin{align}
&\hat{\mathbb{E}}\sup_{0\leq t\leq T}|A^{1,k}_{t}-A^{1,j}_{t}|^{\alpha}\nonumber\\
&\leq C\{\hat{\mathbb{E}}\sup_{0\leq t\leq T}|Y^{1,k}_{s}-Y^{1,j}_{s}|^{\alpha}+\hat{\mathbb{E}}[(\int_{0}^{T}|Z^{1,k}_{s}-Z^{1,j}_{s}|^{2}ds)^{\frac{\alpha}{2}}]\}\rightarrow0,
\end{align}
as $k,j\rightarrow\infty$, which implies that $\{A_{t}^{1,k}\}_{k\geq0}$ is a Cauchy sequence in $S_{G}^{\alpha}(0,T)$.  We denote its limit as $A_{t}^{1}$, i.e. $\hat{\mathbb{E}}[\sup_{0\leq t\leq T}|A^{1,k}_{t}-A^{1}_{t}|^{\alpha}]\rightarrow0$  as $k\rightarrow\infty.$  It is easy to see that $A_{0}^{1}=0$ and $A_{t}^{1}$ is a nondecreasing process, since the sequence $\{A_{t}^{1,k}\}_{k\geq0}$ have the property.

Taking limit on both side of (3.7), then we have
\begin{align}
&Y_{t}^{1}=\xi+\int_{t}^{T}f^{1}(s,w,Y_{s}^{1},Z_{s}^{1})ds-\int_{t}^{T}Z_{s}^{1}dB_{s}+A^{1}_{T}-A^{1}_{t}.
\end{align}
Moreover, $(Y^{1},Z^{1},A^{1})\in \mathcal{S}_{G}^{\alpha}(0,T), \     2\leq\alpha<\beta$.

In the following it remains to prove that $\{-\int_{0}^{t}(Y_{s}^{1}-L_{s})dA_{s}^{1}\}_{t\in [0,T]}$ is a non-increasing $G$-martingale.
Notice that $\{-\int_{0}^{t}(Y_{t}^{1,k}-L_{t})dA_{t}^{1,k}\}_{t\in[0,T]}$ is a non-increasing $G$-martingale. By H\"{o}lder inequality,  (3.10), (3.13) and (3.18), it follows that
\begin{align}
&\hat{\mathbb{E}}[\sup_{0\leq t\leq T}|-\int_{0}^{t}(Y_{t}^{1,k}-L_{t})dA_{t}^{1,k}-(-\int_{0}^{t}(Y_{t}^{1}-L_{t})dA_{t}^{1})|]
\nonumber\\&\leq \hat{\mathbb{E}}[\sup_{0\leq t\leq T}|\int_{0}^{t}(Y_{t}^{1}-Y_{t}^{1,k})dA_{t}^{1,k}|]+ \hat{\mathbb{E}}[\sup_{0\leq t\leq T}|\int_{0}^{t}(Y_{t}^{1}-L_{t})d(A_{t}^{1}-A_{t}^{1,k})|]\nonumber\\&
\leq (\hat{\mathbb{E}}[\sup_{0\leq t\leq T}|Y_{t}^{1}-Y_{t}^{1,k}|^{2}])^{\frac{1}{2}}(\hat{\mathbb{E}}|A_{T}^{1,k}|^{2})^{\frac{1}{2}}+(\hat{\mathbb{E}}[\sup_{0\leq t\leq T}|Y_{t}^{1}-L_{t}|^{2}])^{\frac{1}{2}}(\hat{\mathbb{E}}|A_{T}^{1}-A_{T}^{1,k}|^{2})^{\frac{1}{2}}\rightarrow0,\nonumber
\end{align}
as $k\rightarrow\infty$,  which implies that $\{-\int_{0}^{t}(Y_{t}^{1}-L_{t})dA_{t}^{1}\}_{t\in[0,T]}$ is a non-increasing $G$-martingale.

Finally we obtain a triple $(Y_{t}^{1},Z_{t}^{1},A_{t}^{1})\in \mathcal{S}_{G}^{\alpha}(0,T) $ satisfying the following equation:
\begin{equation}
\begin{cases}
Y_{t}^{1}=\xi+\int_{t}^{T}f^{1}(s,w,Y_{s}^{1},Z_{s}^{1})ds-\int_{t}^{T}Z_{s}^{1}dB_{s}+(A_{T}^{1}-A_{t}^{1}), \\
Y^{1}_{t}\geq L_{t}, \{-\int_{0}^{t}(Y_{s}^{1}-L_{s})dA_{s}^{1}\}_{t\in [0,T]}\      $is$ \      $a$ \      $non-increasing$    \  G-$martingale$.
\end{cases}
\end{equation}
Next, using the same method as that on $X^{0}$, we can construct $X^{1}$ based on $Y^{1}$. Moreover, since $b$ is monotonic on $y$ and $Y^{0}\leq Y^{1}$, we have through the comparison on GSDEs that $X^{0}\leq X^{1}$.

Repeating the same procedure, we get the existence of a sequence $(X^{n},Y^{n},Z^{n},A^{n})$ which is a solution of (3.2) and for any $t\leq T,   \  n\in \mathbb{N}$,
$$X^{n}_{t}\leq X^{n+1}_{t}\leq S_{t},   \    \   L_{t} \leq Y^{n}_{t}\leq Y^{n+1}_{t}\leq U_{t},    q.s..$$
 Moreover,  $X_{t}^{n}\in M_{G}^{2}(0,T), (Y_{t}^{n},Z_{t}^{n}, A_{t}^{n})\in \mathcal{S}_{G}^{\alpha}(0,T)$,  $A^{n}_{t}$ is a nondecreasing  process with $A^{n}_{0}=0$,  where in each step of iteration, we denote $f^{n}(s,w,y,z)=f(s,X_{s}^{n-1},y,z).$

 By the $G$-It\^{o}'s formula to $e^{\lambda t}|Y_{t}^{n}-c|^{2}$, here $c$ is the constant in (H5) and $\lambda>0$ is a constant,  we have
 \begin{align}
&M^{n}_{T}-M^{n}_{t}+e^{\lambda t}|Y^{n}_{t}-c|^{2}+\lambda \int_{t}^{T}e^{\lambda s}|Y_{s}^{n}-c|^{2}ds+\int_{t}^{T}e^{\lambda s}|Z_{s}^{n}|^{2}d\langle B\rangle_{s}\nonumber\\
&\leq e^{\lambda T}|\xi-c|^{2}+2\int_{t}^{T}e^{\lambda s}(Y^{n}_{t}-c)f(s,X^{n}_{s},Y^{n}_{s},Z_{s}^{n})ds，
\end{align}
where $M^{n}_{t}=\int_{0}^{t}e^{\lambda s}(Y^{n}_{t}-c)Z_{s}^{n}dB_{s}-\int_{0}^{t}e^{\lambda s}|Y^{n}_{t}-c|dA_{s}^{n}$.  We claim that $\{M_{t}^{n}\}_{t\in [0,T]}$   is a $G$-martingale. Indeed, note that
$$0\geq-\int_{t}^{T}|Y^{n}_{s}-c|dA_{s}^{n}\geq -\int_{t}^{T}(Y^{n}_{s}-L_{t})dA_{s}^{n}.$$
Thus we can conclude that
$$0\geq\hat{\mathbb{E}}_{t}[-\int_{t}^{T}|Y^{n}_{s}-c|dA_{s}^{n}]\geq \hat{\mathbb{E}}_{t}[-\int_{t}^{T}(Y^{n}_{s}-L_{t})dA_{s}^{n}]=0.$$
It follows that $\{M_{t}^{n}\}_{t\in [0,T]}$   is a $G$-martingale.

Since $f$ satisfies the linear growth condition (H3), we have
 \begin{align}
&2\int_{t}^{T}e^{\lambda s}|Y^{n}_{t}-c|f(s,X^{n}_{s},Y^{n}_{s},Z_{s}^{n})ds
\nonumber\\
&\leq M\int_{t}^{T}e^{\lambda s}(c+1)^{2}ds+(M+\frac{2M^{2}}{\underline{\sigma}^{2}})\int_{t}^{T}e^{\lambda s}|Y_{s}^{n}-c|^{2}ds+\frac{1}{2}\int_{t}^{T}e^{\lambda s}|Z_{s}^{n}|^{2}d\langle B\rangle_{s}.\nonumber
\end{align}

Let $\lambda=M+\frac{2M^{2}}{\underline{\sigma}^{2}}+1$, then
 \begin{align}
&M_{T}^{n}-M_{t}^{n}+e^{\lambda t}|Y^{n}_{t}-c|^{2}+\frac{1}{2}\int_{t}^{T}e^{\lambda s}|Z_{s}^{n}|^{2}d\langle B\rangle_{s}\nonumber\\
&\leq e^{\lambda T}|\xi-c|^{2}+ M\int_{t}^{T}e^{\lambda s}(c+1)^{2}ds，
\end{align}

Taking conditional expectation on both side of (3.22), we have
 \begin{align}
&\hat{\mathbb{E}}_{t}[\int_{t}^{T}e^{\lambda s}|Z_{s}^{n}|^{2}d\langle B\rangle_{s}]\leq 2\hat{\mathbb{E}}_{t}[e^{\lambda T}|\xi-c|^{2}]+2\hat{\mathbb{E}}_{t}[M\int_{t}^{T}e^{\lambda s}(c+1)^{2}ds].
\end{align}

Taking expectation on both side of (3.23),  we have
 \begin{align}
&\underline{\sigma}^{2}\hat{\mathbb{E}}[\int_{0}^{T}e^{\lambda s}|Z_{s}^{n}|^{2}ds]\leq 2\hat{\mathbb{E}}[e^{\lambda T}|\xi-c|^{2}]+2\hat{\mathbb{E}}[M\int_{t}^{T}e^{\lambda s}(c+1)^{2}ds]< C ,
\end{align}
where $C$  is a constant independent of $n$, which implies that $\{Z_{t}^{n}\}_{t\in [0,T]}$ is a bounded process in $M_{G}^{2}(0,T)$  independent of $n$.

Obviously, since $\{X_{t}^{n}\}_{n\in \mathbb{N}}$ and $\{Y_{t}^{n}\}_{n\in \mathbb{N}}$ are increasing and bounded, then there exist two semi-continuous processes $\{X_{t}\}_{t\in[0,T]}$ and $\{Y_{t}\}_{t\in[0,T]}$ such that
$$X_{t}=\lim_{n\rightarrow\infty}X_{t}^{n},    \     \     Y_{t}=\lim_{n\rightarrow\infty}Y_{t}^{n}.$$
Moreover,  we have $|X_{t}^{n}-X_{t}|\downarrow 0$ and $|Y_{t}^{n}-Y_{t}|\downarrow 0$  as $n\rightarrow\infty$.  Notice  that $\{X_{t}^{n}\}_{n\in \mathbb{N}}$ and $\{Y_{t}^{n}\}_{n\in \mathbb{N}}$ are bounded by $S_{t}$ and $U_{t}$  respectively,  which imply that $X_{t},Y_{t}\in S_{G}^{2}(0,T)$.  They also belong to $L_{G}^{2}(\Omega_{T})$,  which is a larger space.    Thanks to the Lemma 2.3, we have
\begin{align}
&\hat{\mathbb{E}}[|X^{n}_{t}-X_{t}|^{2}]\downarrow 0,   \   \ \hat{\mathbb{E}}[|Y^{n}_{t}-Y_{t}|^{2}]\downarrow 0, \end{align}
as $n\rightarrow\infty$. Therefore, by the Lebesgue dominated convergence with respect to $t$, we have
$$\int_{0}^{T}\hat{\mathbb{E}}[|X^{n}_{t}-X_{t}|^{2}]dt\rightarrow 0, \  \   \int_{0}^{T}\hat{\mathbb{E}}[|Y^{n}_{t}-Y_{t}|^{2}]dt\rightarrow 0,$$
as $n\rightarrow\infty$, which imply that $X^{n}_{t}\rightarrow X_{t}, Y^{n}_{t}\rightarrow Y_{t}$ in $M_{G}^{2}(0,T)$.

Now taking the limit on the first equation in (3.2) and using the same method as that on $X_{t}^{0}$, we conclude that $X\in M_{G}^{2}(0,T), Y\in M_{G}^{2}(0,T)$ is a solution to the following forward GSDE:
\begin{align}
&X_{t}=x+\int_{0}^{t}b(s,X_{s},Y_{s})+\int_{0}^{t}\sigma(s,X_{s})dB_{s},
\end{align}
which also implies that $X$ is continuous in $t$.

Let us focus on the backward equation part  of (3.2). For any $n\geq 1$, we have
$$Y^{n}_{t}=\xi+\int_{t}^{T}f(s,X^{n-1}_{s},Y^{n}_{s},Z^{n}_{s})ds-\int_{t}^{T}Z^{n}_{s}dB_{s}+(A^{n}_{T}-A^{n}_{t}).$$
By (3.24),  there exists a real positive constant $C$ independent of $n$ such that $\hat{\mathbb{E}}[\int_{0}^{T}|Z_{t}^{n}|^{2}dt]\leq C$. Therefore the sequence $\{Z^{n}\}_{n\in \mathbb{N}}$ is a Cauchy type in $M_{G}^{2}(0,T)$. Actually,  applying $G$-It\^{o}'s formula to $|Y_{t}^{m}-Y_{t}^{n}|^{2}$, we can obtain
\begin{align}
&|Y_{t}^{m}-Y_{t}^{n}|^{2}+\int_{t}^{T}|Z_{s}^{m}-Z_{s}^{n}|^{2}d\langle B\rangle_{s}\nonumber\\
&=2\int_{t}^{T}(Y_{s}^{m}-Y_{s}^{n})[f(s,X_{s}^{m-1},Y_{s}^{m},Z_{s}^{m})-f(s,X_{s}^{n-1},Y_{s}^{n},Z_{s}^{n})]ds\nonumber\\
&+2\int_{t}^{T}(Y_{s}^{m}-Y_{s}^{n})d(A_{s}^{m}-A_{s}^{n})-2\int_{t}^{T}(Y_{t}^{m}-Y_{t}^{n})(Z_{s}^{m}-Z_{s}^{n})dB_{s}\nonumber\\
&\leq 2\int_{t}^{T}(Y_{s}^{m}-Y_{s}^{n})[f(s,X_{s}^{m-1},Y_{s}^{m},Z_{s}^{m})-f(s,X_{s}^{n-1},Y_{s}^{n},Z_{s}^{n})]ds\nonumber\\
&+2\int_{t}^{T}(Y_{s}^{m}-Y_{s}^{n})^{+}dA_{t}^{m}+2\int_{t}^{T}(Y_{s}^{m}-Y_{s}^{n})^{-}dA_{s}^{n}-2\int_{t}^{T}(Y_{s}^{m}-Y_{s}^{n})(Z_{s}^{m}-Z_{s}^{n})dB_{s}.
\end{align}
Let $M_{t}^{m,n}=2\int_{0}^{t}(Y_{s}^{m}-Y_{s}^{n})(Z_{s}^{m}-Z_{s}^{n})dB_{s}-2\int_{0}^{t}(Y_{s}^{m}-Y_{s}^{n})^{+}dA_{t}^{m}-2\int_{0}^{t}(Y_{s}^{m}-Y_{s}^{n})^{-}dA_{s}^{n},$ which is a $G$-martingale.  Then we have
\begin{align}
&M_{T}^{m,n}-M_{t}^{m,n}+|Y_{t}^{m}-Y_{t}^{n}|^{2}+\int_{t}^{T}|Z_{s}^{m}-Z_{s}^{n}|^{2}d\langle B\rangle_{s}\nonumber\\
&\leq 2\int_{t}^{T}(Y_{s}^{m}-Y_{s}^{n})[f(s,X_{s}^{m-1},Y_{s}^{m},Z_{s}^{m})-f(s,X_{s}^{n-1},Y_{s}^{n},Z_{s}^{n})]ds.\nonumber
\end{align}

Taking expectation we may conclude
\begin{align}
&\underline{\sigma}^{2}\hat{\mathbb{E}}[\int_{0}^{T}|Z_{s}^{m}-Z_{s}^{n}|^{2}ds]\leq \hat{\mathbb{E}}[\int_{0}^{T}|Z_{s}^{m}-Z_{s}^{n}|^{2}d\langle B\rangle_{s}]\nonumber\\
&\leq(\int_{0}^{T}\hat{\mathbb{E}}|Y_{s}^{m}-Y_{s}^{n}|^{2}ds)^{\frac{1}{2}}(\int_{0}^{T}\hat{\mathbb{E}}|f(s,X_{s}^{m-1},Y_{s}^{m},Z_{s}^{m})-f(s,X_{s}^{n-1},Y_{s}^{n},Z_{s}^{n})|^{2}ds)^{\frac{1}{2}},
\end{align}
Notice that  $|f(t,x,y,z)|\leq K(1+|y|+|z|)$, $\{Z^{n}\}_{n\in \mathbb{N}}$ is bounded in $M_{G}^{2}(0,T)$, $\{X^{n}\}_{n\in \mathbb{N}}$ and $\{Y^{n}\}_{n\in \mathbb{N}}$ are bounded by the process $S$  and $ U$  respectively,  then the second factor in above inequality is bounded by constant independent of $m$ and $n$. Therefore, $\{Z^{n}\}_{n\in \mathbb{N}}$ is a Cauchy sequence in $M_{G}^{2}(0,T)$ through the convergence of $\{Y^{n}\}_{n\in \mathbb{N}}$ to $Y$ in $M_{G}^{2}(0,T)$. So let us set $Z_{t}:=\lim_{n\rightarrow\infty}Z^{n}_{t}$ in $M_{G}^{2}(0,T)$.

By  (3.2), we have
$$A_{t}^{n}=Y_{0}^{n}-Y_{t}^{n}-\int_{0}^{t}f(s,X_{s}^{n-1},Y_{s}^{n}, Z_{s}^{n})ds+\int_{0}^{t}Z_{s}^{n}dB_{s}.$$
Thus
\begin{align}
&\hat{\mathbb{E}}[|A_{t}^{m}-A_{t}^{n}|^{2}]\leq 4\{\hat{\mathbb{E}}[|Y_{0}^{m}-Y_{0}^{n}|^{2}]+|Y_{t}^{m}-Y_{t}^{n}|^{2}]\nonumber\\
&+T\int_{0}^{T}\hat{\mathbb{E}}|f(s,X_{s}^{m-1},Y_{s}^{m},Z_{s}^{m})-f(s,X_{s}^{n-1},Y_{s}^{n},Z_{s}^{n})|^{2}ds+\hat{\mathbb{E}}[|\int_{0}^{T}(Z_{t}^{m}-Z_{t}^{n})dB_{t}|^{2}]\}
\nonumber\\
&=:4(I_{1}+I_{2}+I_{3}).
\end{align}

By (3.25), we have
\begin{align}
&I_{1}\leq\hat{\mathbb{E}}[|Y_{0}^{m}-Y_{0}|^{2}]+\hat{\mathbb{E}}[|Y_{0}^{n}-Y_{0}|^{2}]+\hat{\mathbb{E}}[|Y_{t}^{m}-Y_{t}|^{2}]+\hat{\mathbb{E}}[|Y_{t}^{n}-Y_{t}|^{2}]\rightarrow 0,
\end{align}
as $m,n\rightarrow\infty$.

By Lemma 2.1 and $\{Z^{n}\}_{n\in \mathbb{N}}$ is a Cauchy sequence in $M_{G}^{2}(0,T)$, we have
\begin{align}
&I_{3}\leq \bar{\sigma}^{2}\hat{\mathbb{E}}[\int_{0}^{T}|Z_{t}^{m}-Z_{t}^{n}|^{2}ds]\rightarrow 0,
\end{align}
as $m,n\rightarrow\infty$.

By Lemma 2.5 and$\{X^{n}\}_{n\in \mathbb{N}}$,  $\{Y^{n}\}_{n\in \mathbb{N}}$ and $\{Z^{n}\}_{n\in \mathbb{N}}$ are convergence in $M_{G}^{2}(0,T)$
\begin{align}
&I_{2}\leq 3T\int_{0}^{T}\hat{\mathbb{E}}[|f(s,X_{s}^{m-1},Y_{s}^{m},Z_{s}^{m})-f_{k}(s,X_{s}^{m-1},Y_{s}^{m},Z_{s}^{m})|^{2}]ds\nonumber\\
&+\int_{0}^{T}\hat{\mathbb{E}}[|f_{k}(s,X_{s}^{m-1},Y_{s}^{m},Z_{s}^{m})-f_{k}(s,X_{s}^{n-1},Y_{s}^{n},Z_{s}^{n})|^{2}]ds\nonumber\\
&+\int_{0}^{T}\hat{\mathbb{E}}[|f_{k}(s,X_{s}^{n-1},Y_{s}^{n},Z_{s}^{n})-f(s,X_{s}^{n-1},Y_{s}^{n},Z_{s}^{n})|^{2}]ds\rightarrow0,
\end{align}
as $k,m,n\rightarrow\infty$,  where $f_{k}$ is the sequence defined in Lemma 2.5.

Combining  (3.29)-(3.32),  then  $\{A^{n}\}_{n\in \mathbb{N}}$ is a Cauchy sequence in $L_{G}^{2}(\Omega_{T})$.  Let us set $A_{t}:=\lim_{n\rightarrow\infty}A_{t}^{n}$ in $L_{G}^{2}(\Omega_{T})$.  From this, we can show that $\{Y^{n}\}_{n\in \mathbb{N}}$ is a Cauchy sequence in $S_{G}^{2}(\Omega_{T})$. Actually, by (3.27), we have
\begin{align}
&\hat{\mathbb{E}}[\sup_{0\leq t \leq T}|Y_{t}^{m}-Y_{t}^{n}|^{2}]\nonumber\\
&=2\hat{\mathbb{E}}[\int_{0}^{T}(Y_{s}^{m}-Y_{s}^{n})(f(s,X_{s}^{m-1},Y_{s}^{m},Z_{s}^{m})-f(s,X_{s}^{n-1},Y_{s}^{n},Z_{s}^{n}))ds]\nonumber\\
&+2\hat{\mathbb{E}}[\int_{0}^{T}(Y_{s}^{m}-Y_{s}^{n})d(A_{s}^{m}-A_{s}^{n})]+2\hat{\mathbb{E}}[\int_{0}^{T}(Y_{t}^{m}-Y_{t}^{n})(Z_{s}^{m}-Z_{s}^{n})dB_{s}]\nonumber\\
&\leq 2(\hat{\mathbb{E}}[\int_{0}^{T}|Y_{t}^{m}-Y_{t}^{n}|^{2}ds])^{\frac{1}{2}}(\hat{\mathbb{E}}[\int_{0}^{T}|f(s,X_{s}^{m-1},Y_{s}^{m},Z_{s}^{m})-f(s,X_{s}^{n-1},Y_{s}^{n},Z_{s}^{n})|^{2}ds])^{\frac{1}{2}}\nonumber\\
&+2\hat{\mathbb{E}}[\sup_{0\leq t \leq T}|Y_{t}^{m}-Y_{t}^{n}|\cdot|A_{T}^{m}-A_{T}^{n}|]+2\hat{\mathbb{E}}[(\int_{0}^{T}|Y_{s}^{m}-Y_{s}^{n}|^{2}|Z_{s}^{m}-Z_{s}^{n}|^{2}ds)^{\frac{1}{2}}]\nonumber\\
&\leq 2(\hat{\mathbb{E}}[\int_{0}^{T}|Y_{t}^{m}-Y_{t}^{n}|^{2}ds])^{\frac{1}{2}}(\hat{\mathbb{E}}[\int_{0}^{T}|f(s,X_{s}^{m-1},Y_{s}^{m},Z_{s}^{m})-f(s,X_{s}^{n-1},Y_{s}^{n},Z_{s}^{n})|^{2}ds])^{\frac{1}{2}}\nonumber\\
&+\varepsilon \hat{\mathbb{E}}[\sup_{0\leq t \leq T}|Y_{t}^{m}-Y_{t}^{n}|^{2}]+\frac{1}{\varepsilon}\hat{\mathbb{E}}[|A_{T}^{m}-A_{T}^{n}|^{2}]+\varepsilon\hat{\mathbb{E}}[\sup_{0\leq t \leq T}|Y_{t}^{m}-Y_{t}^{n}|^{2}]+\frac{1}{\varepsilon}\hat{\mathbb{E}}[\int_{0}^{T}|Z_{s}^{m}-Z_{s}^{n}|^{2}ds],\nonumber
\end{align}
where in the above  inequalities we have used the H\"{o}lder inequality, Lemma 2.1 and Young inequality. Choosing $\varepsilon$ small enough, then we have
\begin{align}
&\hat{\mathbb{E}}\sup_{0\leq t \leq T}|Y_{t}^{m}-Y_{t}^{n}|^{2}\nonumber\\
&\leq C(\hat{\mathbb{E}}[\int_{0}^{T}|Y_{t}^{m}-Y_{t}^{n}|^{2}ds)^{\frac{1}{2}}]+\frac{1}{\varepsilon}\hat{\mathbb{E}}[|A_{T}^{m}-A_{T}^{n}|^{2}]+\frac{1}{\varepsilon}\hat{\mathbb{E}}[\int_{0}^{T}|Z_{s}^{m}-Z_{s}^{n}|^{2}ds]\rightarrow 0,
\end{align}
as $m,n\rightarrow\infty$, which implies that $Y_{t}^{n}\rightarrow Y_{t}$ in $S_{G}^{2}(0,T)$ as $n\rightarrow\infty.$  Going back to (3.29), it is easy to see that $A_{t}^{n}\rightarrow A_{t}$ in $S_{G}^{2}(0,T)$ as $n\rightarrow\infty.$

By (H3), Lemma 2.5 and the convergence of $\{X^{n}\}_{n\in \mathbb{N}},  \{Y^{n}\}_{n\in \mathbb{N}}, \{Z^{n}\}_{n\in \mathbb{N}}$, similar to that in (3.32), we have
\begin{align}
&\hat{\mathbb{E}}[\int_{0}^{T}|f(s,X_{s}^{n-1},Y_{s}^{n},Z_{s}^{n})-f(s,X_{s},Y_{s},Z_{s})|^{2}ds]\rightarrow 0
\end{align}
Taking limit on both side of (3.2), then we have
\begin{align}
&Y_{t}=\xi+\int_{t}^{T}f(s,X_{s},Y_{s}, Z_{s})ds-\int_{t}^{T}Z_{s}dB_{s}+A_{T}-A_{t}    \ \  q.s..
\end{align}
Together with  (3.26),  then the process $(X,Y,Z,A )$ satisfy (3.1). Moreover, $X_{t}\in M_{G}^{2}(0,T)$,  $(Y_{t},Z_{t}, A_{t})\in \mathcal{S}_{G}^{2}(0,T).$

In the following, it remains to prove that $\{-\int_{0}^{t}(Y_{s}-L_{s})dA_{s}\}_{t\in[0,T]}$  is a non-increasing $G$-martingale. Note that $\{-\int_{0}^{t}(Y_{t}^{n}-L_{t})dA_{t}^{n}\}_{t\in[0,T]}$ is a non-increasing  $G$-martingale. Moreover,
\begin{align}
&\hat{\mathbb{E}}[\sup_{0\leq t\leq T}|-\int_{0}^{t}(Y_{t}^{n}-L_{t})dA_{t}^{n}-(-\int_{0}^{t}(Y_{t}-L_{t})dA_{t})|]
\nonumber\\&\leq \hat{\mathbb{E}}[\sup_{0\leq t\leq T}|\int_{0}^{t}(Y_{t}-Y_{t}^{n})dA_{t}^{n}|]+ \hat{\mathbb{E}}[\sup_{0\leq t\leq T}|\int_{0}^{t}(Y_{t}-L_{t})d(A_{t}-A_{t}^{n})|]\nonumber\\&
\leq (\hat{\mathbb{E}}[\sup_{0\leq t\leq T}|Y_{t}-Y_{t}^{n}|^{2}])^{\frac{1}{2}}](\hat{\mathbb{E}}[|A_{T}^{n}|^{2}])^{\frac{1}{2}}+(\hat{\mathbb{E}}[\sup_{0\leq t\leq T}|Y_{t}-L_{t}|^{2}])^{\frac{1}{2}}(\hat{\mathbb{E}}[|A_{T}-A_{T}^{n}|^{2}])^{\frac{1}{2}}\rightarrow0,
\end{align}
as $n\rightarrow\infty$,  which implies that $\{-\int_{0}^{t}(Y_{t}-L_{t})dA_{t}\}_{t\in[0,T]}$ is a non-increasing $G$-martingale.

Furthermore, we have the following result.
\begin{theorem}
Suppose that $\xi, b, h, f,g,  \sigma$ satisfy (H1)-(H4), $L$ satisfies (H5). Then the RFBGSDE (1.1) has at least one solution $(X,Y,Z,A)$.
\end{theorem}

\begin{remark}
The same proof works also when the terminal condition $\xi$ is replaced by $\Phi(\xi)$, where $\Phi$ is a continuous bounded increasing function.
\end{remark}

\begin{remark}
In our assumptions, $f$ has sublinear  growth independent of $x$. If we assume that
$$|f(t,x,y,z)|\leq M(1+|x|+|y|+|z|),  \forall t\in [0,T]\    \   x,y,z\in R,$$
then $b$ should have a sunlinear growth independent of $y$, i.e.
$$b(s,x,y)\leq M(1+|x|).$$
We can also construct a sequence of Lipschitz-continuous functions to approximate them. However, in this case, we should first construct a solution of some forward GSDEs, which is different from that in this paper.
\end{remark}


\renewcommand{\theequation}{A\arabic{equation}}
\setcounter{equation}{0}
%
%

\end{document}